\documentclass[12pt, a4]{amsart}
\usepackage{verbatim}
\numberwithin{equation}{section}
\newtheorem{thm}{Theorem}

\newcommand\cal{\mathcal}%%%%                                                                            
\begin{document}
\title
{On random multiple Dirichlet series}
\author{Gautami Bhowmik and Kohji Matsumoto}
\date{}
\maketitle
%%%%%%%%%%%%%%%%%%%%%%%%%%%%%%%%%%%%%%%%%%%%%%%%%%%%%%%%%%%%%%%%%%%%%%%%%%%%%%  
%                                                                               
%   1                                                                           
%                                                                               
%%%%%%%%%%%%%%%%%%%%%%%%%%%%%%%%%%%%%%%%%%%%%%%%%%%%%%%%%%%%%%%%%%%%%%%%%%%%%%  
\section{Introduction}

A sequences of independent random variables, each of which takes only the
values $\pm 1$, is called a Rademacher sequence.
Dirichlet series or Taylor series whose coefficients include Rademacher type of 
variables are studied, for example, in Kahane \cite{Kah}.

In the present paper we study a multiple analogue of Rademacher sequences and
associated Dirichlet series.  Let
$m_1,\ldots,m_r\in{\Bbb N}$, and let $\varepsilon_{m_1,\ldots, m_r}$ be independent
random variables which take only the values $\pm 1$ and are
defined on a certain probability space $(\Omega,P)$. Let 

\begin{align}\label{1-1}
P(\omega\in\Omega\;|\;\varepsilon_{m_1,\ldots, m_r}(\omega)=1)
=P(\omega\in\Omega\;|\;\varepsilon_{m_1,\ldots, m_r}(\omega)=-1)=\frac{1}{2}.
\end{align}
The sequence $\{\varepsilon_{m_1,\ldots, m_r}\;|\;m_1,\ldots,m_r\in{\Bbb N}\}$
may be called a Rademacher sequence of $r$-tuple indices.
The associated multiple Dirichlet series can be defined as
\begin{align}\label{1-2}
F_r(s,\omega)=\sum_{m_1=1}^{\infty}\cdots\sum_{m_r=1}^{\infty}
\frac{\varepsilon_{m_1,\ldots, m_r}(\omega)X(m_1,\ldots,m_r)}{(m_1+\cdots+m_r)^s},
\end{align}
where $s\in{\Bbb C}$ and $X(m_1,\ldots,m_r)\in{\Bbb R}$.

As an example, we will treat the case where
$X(m_1,\ldots,m_r)$ is a product of the von Mangoldt function $\Lambda(\cdot)$,
that is
\begin{align}\label{1-3}                                                              
\Phi_r(s,\omega)=\sum_{m_1=1}^{\infty}\cdots\sum_{m_r=1}^{\infty}                         
\frac{\varepsilon_{m_1,\ldots, m_r}(\omega)\Lambda(m_1)\cdots\Lambda(m_r)}
{(m_1+\cdots+m_r)^s}.      
\end{align}
The Dirichlet series
\begin{align}\label{1-4}
\Phi_r(s)=\sum_{m_1=1}^{\infty}\cdots\sum_{m_r=1}^{\infty}                        
\frac{\Lambda(m_1)\cdots\Lambda(m_r)}{(m_1+\cdots+m_r)^s}       
\end{align}
is connected with the problem of expressing positive integers as a sum of $r$
prime numbers and the analytic properties of $\Phi_r(s)$ have been
studied actively in the recent past (see \cite{Bh} for references).

In the present paper we discuss the natural boundary of \eqref{1-2} and 
\eqref{1-3}, and our
result for the random series could be compared with those of the correspondind
ordinary  Dirichlet series.

%%%%%%%%%%%%%%%%%%%%%%%%%%%%%%%%%%%%%%%%%%%%%%%%%%%%%%%%%%%%%%%%%%%%%%%%%%%%%
%
%   2
%
%%%%%%%%%%%%%%%%%%%%%%%%%%%%%%%%%%%%%%%%%%%%%%%%%%%%%%%%%%%%%%%%%%%%%%%%%%%%%
\section{The natural boundary}

First we write $F_r(s,\omega)$ in a form of single Dirichlet series:
\begin{align}\label{2-1}
F_r(s,\omega)=\sum_{n=1}^{\infty}Y_r(n,\omega)n^{-s},
\end{align}
where
\begin{align}\label{2-2}
Y_r(n,\omega)=\sum_{m_1+\cdots+m_r=n}\varepsilon_{m_1,\ldots, m_r}(\omega)
    X(m_1,\ldots,m_r).
\end{align}
Since $\epsilon_{m_1,\ldots, m_r}$ are independent variables, we see that
$Y_r(n,\omega)$ are also independent of each other.

For each fixed $\omega\in\Omega$, the series \eqref{2-1} is a Dirichlet series,
so we can find its abscissa of convergence $\sigma_c(F_r,\omega)$.   Then, according to
the zero-one law, we can find a constant $\sigma_c(F_r)$ such that
$\sigma_c(F_r,\omega)=\sigma_c(F_r)$ almost surely 
(see \cite[Section 6 of Chapter 4]{Kah}).
In general $-\infty\leq\sigma_c(F_r) \leq\infty$, but here we assume
$-\infty<\sigma_c(F_r)<\infty$.   Our main result in this section is as follows.

\begin{thm}\label{th-1}
The line $\Re s=\sigma_c(F_r)$ is the natural boundary of $F_r(s,\omega)$ almost
surely.
\end{thm}

\begin{proof}
We say that $Y_r(n,\omega)$ is symmetric if $Y_r(n,\omega)$ and $-Y_r(n,\omega)$
have the same distribution.   We prove that $Y_r(n,\omega)$ is symmetric;
and the conclusion then follows immediately from the expression
\eqref{2-1} and \cite[Theorem 4 in Section 6 of Chapter 4]{Kah}.

Let $L(n)$ be the number of tuples $(m_1,\ldots, m_r)$ satisfying
$m_1+\cdots+m_r=n$.   
Let ${\bf b}$ be any tuple of $L(n)$ elements,  each element being  1 or
$-1$.   Write ${\bf b}=(b_{m_1,\ldots, m_r})_{m_1+\cdots+m_r=n}$, where
$b_{m_1,\ldots, m_r}\in\{\pm 1\}$.
Define
$$
Z_r(n,{\bf b})=\sum_{m_1+\cdots+m_r=n}b_{m_1,\ldots, m_r}X(m_1,\ldots, m_r).
$$
Let $A$ be any Borel subset of ${\Bbb R}$, and let ${\cal B}(A)$ be
the set of all ${\bf b}$ such that $Z_r(n,{\bf b})\in A$.   Then
\begin{align}\label{2-3}
\begin{split}
&P(\omega\in\Omega\;|\; Y_r(n,\omega)\in A)\\
&=\sum_{{\bf b}\in{\cal B}(A)}P\left(\omega\in\Omega\;|\;
  (\varepsilon_{m_1,\ldots, m_r}(\omega))_{m_1+\cdots+m_r=n}={\bf b}\right).
\end{split}
\end{align}
Since $\varepsilon_{m_1,\ldots, m_r}$ are independent, we have
\begin{align}\label{2-4}
P\left(\omega\in\Omega\;|\;                                 
  (\varepsilon_{m_1,\ldots, m_r}(\omega))_{m_1+\cdots+m_r=n}={\bf b}\right)
=2^{-L(n)},
\end{align}
and the same equality holds if we replace ${\bf b}$ by $-{\bf b}$.
Therefore the right-hand side of \eqref{2-3} is equal to
\begin{align*}
&\sum_{{\bf b}\in{\cal B}(A)}P\left(\omega\in\Omega\;|\;                                 
  (\varepsilon_{m_1,\ldots, m_r}(\omega))_{m_1+\cdots+m_r=n}=-{\bf b}\right)\\
&=P(\omega\in\Omega\;|\; Y_r(n,\omega)\in -A).
\end{align*}
This implies that two random variables $Y_r(n,\omega)$ and $-Y_r(n,\omega)$
have the same distribution.   
\end{proof}

%%%%%%%%%%%%%%%%%%%%%%%%%%%%%%%%%%%%%%%%%%%%%%%%%%%%%%%%%%%%%%%%%%%%%%%%%%%%%
%
%  3
%
%%%%%%%%%%%%%%%%%%%%%%%%%%%%%%%%%%%%%%%%%%%%%%%%%%%%%%%%%%%%%%%%%%%%%%%%%%%%
\section{The case of Goldbach generating functions}

In this section we consider the special case \eqref{1-3}.  The function $\Phi_r(s)$ defined by \eqref{1-4}
is the generating function of
\begin{align}\label{3-1}
G_r(n)=\sum_{m_1+\cdots+m_r=n}\Lambda(m_1)\cdots\Lambda(m_r).
\end{align} 
%\begin{comment}
% We first mention the
%history of related theory.   The function $\Phi_r(s)$ defined by \eqref{1-4}
%is the generating function of
%\begin{align}\label{3-1}
%G_r(n)=\sum_{m_1+\cdots+m_r=n}\Lambda(m_1)\cdots\Lambda(m_r).
%\end{align}
%In particular, when $r=2$, the function $G_2(n)$ is clearly related with the famous
%Goldbach problem, and has been studied from the days of Hardy and Littlewood.
%Under the assumption of the Riemann Hypothesis (RH) for the Riemann
%zeta-function $\zeta(s)$, Fujii \cite{Fuj91a} proved the formula
%\begin{align}\label{3-2}
%\sum_{n\leq X}G_2(n)=\frac{1}{2}X^2+O(X^{3/2}),
%\end{align}
%and the above error term has been further refined by Fujii \cite{Fuj91b} and
%Bhowmik and Schlage-Puchta \cite{BSP2}.   In \cite{BSP2} it is shown 
%(under RH) that
%\begin{align}\label{3-3}
%\sum_{n\leq X}G_2(n)=\frac{1}{2}X^2-H(X)+E_2(X)
%\end{align}
%for any large positive $X$, where
%$$
%H(X)=2\sum_{\rho}\frac{X^{1+\rho}}{\rho(1+\rho)}
%$$
%($\rho$ runs over all non-trivial zeros of $\zeta(s)$), and $E_2(X)$ is the 
%error term satisfying 
%\begin{align}\label{3-4}
%E_2(X)=O(X(\log X)^5), \qquad E_2(X)=\Omega(x\log\log X).
%\end{align}
%This confirms Conjecture 2.2 of \cite{EM} in a more precise form. 
%\end{comment}
The series $\Phi_r(s)$ was introduced in \cite{EM} where it was
conjectured that:

\bigskip

\noindent (C-$r$) the line $\Re s=r-1$ would be the natural boundary of 
$\Phi_r(s)$.
\bigskip
%\begin{comment}
%\noindent (C-2) The line $\Re s=1$ would be the natural boundary of $\Phi_2(s)$.
%\bigskip
%\end{comment} 
   
The conjecture is out of reach for the moment. However one can say
more under the RH and other
reasonable
ones on the zeros of the Riemann zeta function like the following due to 
Fujii \cite{Fuj91a} :

\bigskip

\noindent (Z) Let ${\cal I}$ be the set of all imaginary parts of 
non-trivial zeros of $\zeta(s)$. If $\gamma_j\in{\cal I}$
($1\leq j\leq 4$) and $\gamma_1+\gamma_2=\gamma_3+\gamma_4\neq 0$, then
$\{\gamma_1,\gamma_2\}=\{\gamma_3,\gamma_4\}$. 
\bigskip

\noindent In \cite{EM} and  \cite{BSP3}, it was
shown that (C-2) is indeed true under the RH and 
a certain quantitative version of (Z)
or (Z) itself. The case of (C-r) is shown to be true if and only if
(C-2) is. 

%\begin{comment}
%
%\noindent This conjecture would follow from (C-2) under RH, as was shown in
%\cite{BSP3}.   Also in \cite{BSP3} an asymptotic formula for the sum of
%$G_r(n)$ with an $\Omega$-result on the error term was given.
%\end{comment}
The reason of introducing the random series \eqref{1-3} is to observe
the situation from a stochastic viewpoint.   Let
\begin{align}\label{3-5}                                                               
G_r(n,\omega)=\sum_{m_1+\cdots+m_r=n}\varepsilon_{m_1,\ldots,m_r}(\omega)
\Lambda(m_1)\cdots\Lambda(m_r).                          
\end{align}
Then 
\begin{align}\label{3-6}
|G_r(n,\omega)|\leq n^{r-1}(\log n)^r,
\end{align}
and hence \eqref{1-3}, which can be
written as
\begin{align}\label{3-7}
\Phi_r(s,\omega)=\sum_{n=1}^{\infty}G_r(n,\omega)n^{-s},
\end{align}
is absolutely convergent for $\Re s>r$.   We will show:

\begin{thm}\label{th-2}
The random series \eqref{3-7} converges for $\Re s>r/2$ almost surely.
\end{thm}

We will prove this theorem in the next section.   
In particular this theorem implies $\sigma_c(\Phi_r,\omega)\leq r/2$ almost surely,
and hence $\sigma_c(\Phi_r)\leq r/2$.   The lower-bound of $\sigma_c(\Phi_r)$
will be discussed in Section 5.

%%%%%%%%%%%%%%%%%%%%%%%%%%%%%%%%%%%%%%%%%%%%%%%%%%%%%%%%%%%%%%%%%%%%%%%%%%%%
%
%    4
%
%%%%%%%%%%%%%%%%%%%%%%%%%%%%%%%%%%%%%%%%%%%%%%%%%%%%%%%%%%%%%%%%%%%%%%%%%%%%
\section {Proof of Theorem \ref{th-2}}

Let $E(\cdot)$ denote the expected value, and $V(\cdot)$ the variance.
Write $X_n=G_r(n,\omega)n^{-s}$.   From \eqref{3-6} we see that
\begin{align}\label{4-1}
\int_{\Omega}|X_n|^2 d\omega\leq \left(\frac{n^{r-1}(\log n)^r}{n^{\sigma}}
\right)^2\int_{\Omega} d\omega<+\infty,
\end{align}
that is $X_n\in L^2(\Omega)$.

We use Theorem 2 in \cite[Section 2 of Chapter 3]{Kah}, which asserts that if
$X_n\in L^2(\Omega)$ are independent random variables satisfying
$E(X_n)=0$ and $\sum_{n=1}^{\infty}V(X_n)<+\infty$, then $\sum_{n=1}^{\infty}X_n$
converges almost surely.

Since
\begin{align}\label{4-2}
E(X_n)=\frac{1}{n^s}\sum_{m_1+\cdots+m_r=n}\Lambda(m_1)\cdots\Lambda(m_r)
\int_{\Omega}\varepsilon_{m_1,\ldots,m_r}(\omega)d\omega,
\end{align}
and $\varepsilon_{m_1,\ldots,m_r}$ is an element of a Rademacher sequence, it
is obvious that $E(X_n)=0$ for any $n$.   
Next consider 
$V(X_n)=E(|X_n-E(X_n)|^2)=E(|X_n|^2)$.   We see that
\begin{align}\label{4-3}
\begin{split}
V(X_n)&=\frac{1}{n^{2\sigma}}\sum_{{m_1+\cdots+m_r=n}\atop{m'_1+\cdots+m'_r=n}}
\Lambda(m_1)\cdots\Lambda(m_r)\Lambda(m'_1)\cdots\Lambda(m'_r)\\
&\times\int_{\Omega}\varepsilon_{m_1,\ldots,m_r}(\omega)
\varepsilon_{m'_1,\ldots,m'_r}(\omega)d\omega.
\end{split}
\end{align}
If $(m_1,\ldots,m_r)\neq(m'_1,\ldots,m'_r)$, then
$$
P(\omega\in\Omega\;|\;\varepsilon_{m_1,\ldots,m_r}(\omega)=\pm 1,
\;\varepsilon_{m'_1,\ldots,m'_r}(\omega)=\pm 1)=\frac{1}{4}
$$
for any choice of double signs, hence
$$
\int_{\Omega}\varepsilon_{m_1,\ldots,m_r}(\omega)                                  
\varepsilon_{m'_1,\ldots,m'_r}(\omega)d\omega=0.
$$
Therefore
\begin{align}\label{4-4}
\begin{split}
V(X_n)&=\frac{1}{n^{2\sigma}}\sum_{m_1+\cdots+m_r=n}\Lambda(m_1)^2\cdots
\Lambda(m_r)^2 \int_{\Omega}\varepsilon_{m_1,\ldots,m_r}(\omega)^2d\omega\\
&=\frac{1}{n^{2\sigma}}\sum_{m_1+\cdots+m_r=n}\Lambda(m_1)^2\cdots
\Lambda(m_r)^2,
\end{split}
\end{align}
and the sum on the right-hand side can be estimated as
$\leq n^{r-1}(\log n)^{2r}$.   Therefore
\begin{align}\label{4-5}
\sum_{n=1}^{\infty}V(X_n)\leq \sum_{n=1}^{\infty}n^{r-1-2\sigma}(\log n)^{2r},
\end{align}
which is convergent if $\sigma> r/2$.   This implies the assertion of 
Theorem \ref{th-2}.

%%%%%%%%%%%%%%%%%%%%%%%%%%%%%%%%%%%%%%%%%%%%%%%%%%%%%%%%%%%%%%%%%%%%%%%%%%%%%%
%
%  5
%
%%%%%%%%%%%%%%%%%%%%%%%%%%%%%%%%%%%%%%%%%%%%%%%%%%%%%%%%%%%%%%%%%%%%%%%%%%%%%
\section{The  abscissa of  convergence of $\Phi_r(s,\omega)$}

Theorem \ref{th-2} gives the upper-bound $\sigma_c(\Phi_r)\leq r/2$.
In this section we comment on its lower-bound.   Since $G_r$ is symmetric
(as we have seen in the proof of Theorem \ref{th-1}), we have the criterion that
$\sum_{n=1}^{\infty}X_n$ converges almost surely if and only if 
$\sum_{n=1}^{\infty}V(X'_n)$
converges, where
\begin{align}\label{5-1}
X'_n(\omega)=\left\{
   \begin{array}{lll}
     X_n(\omega) & {\rm if} & |X_n(\omega)|\leq 1 \\
     X_n(\omega)/|X_n(\omega)| & {\rm if} & |X_n(\omega)|>1 
\end{array}
\right.
\end{align}
(\cite[Theorem 7 in Section 5 of Chapter 3]{Kah}).

In the case $r=1$, we have $X_n(\omega)=\varepsilon_n(\omega)\Lambda(n)n^{-s}$.
Hence if $\sigma>0$ then $|X_n(\omega)|\leq 1$ for sufficiently large $n$, so
we may assume $X'_n=X_n$ if $\sigma>0$.   From \eqref{4-4} we have
\begin{align}\label{5-2}
\sum_{n=1}^{\infty}V(X_n)=\sum_{n=1}^{\infty}\frac{\Lambda(n)^2}{n^{2\sigma}},
\end{align}
which is convergent if and only if $\sigma>1/2$.   This implies that
$\sigma_c(\Phi_1)=1/2$, and hence, by Theorem \ref{th-1}, the line $\Re s=1/2$ is 
the natural boundary of $\Phi_1(s,\omega)$ almost surely.

Next consider the case $r=2$.   Then
$$
X_n(\omega)=\frac{1}{n^s}\sum_{m_1+m_2=n}\varepsilon_{m_1,m_2}\Lambda(m_1)
\Lambda(m_2).
$$
If $\sigma>1$ then we may assume $X'_n=X_n$ for sufficiently large $n$.
>From \eqref{4-4} we have
\begin{align}\label{5-3}                                                              
\sum_{n=1}^{\infty}V(X_n)=\sum_{m_1+m_2=n}\frac{\Lambda(m_1)^2\Lambda(m_2)^2}
{n^{2\sigma}},        
\end{align}
whose right-hand side is
\begin{align}\label{5-4}
\geq \frac{1}{n^{2\sigma}}\sum_{{p_i,p_2:{\rm prime}}\atop{p_1+p_2=n}}
(\log p_1)^2(\log p_2)^2
\gg \frac{(\log n)^2}{n^{2\sigma}}\sum_{{p_i,p_2:{\rm prime}}\atop{p_1+p_2=n}}1,
\end{align}
because one of $p_1$ or $p_2$ is $\geq n/2$.
If the Hardy-Littlewood conjectural asymptotic formula for the Goldbach
conjecture is true, then the last sum of \eqref{5-4} is
$\gg n(\log n)^{-2}$.   This implies that \eqref{5-3} diverges for $\sigma=1$.
Though $\sigma=1$ is out of the range where $X'_n=X_n$ is valid, this argument
suggests that $\sigma_c(\Phi_2)=1$, and hence $\Re s=1$ would be the natural boundary 
of $\Phi_2(s,\omega)$ almost surely.

The above observation in the cases $r=1$ and $r=2$ further suggests that, even
in the case $r\geq 3$, perhaps $\sigma_c(\Phi_r)=r/2$ would hold.   
It is to be noted that Theorem \ref{th-2} says,
in the case $r\geq 3$, $\Phi_r(s,\omega)$ is convergent almost surely beyond the
(conjectural) natural boundary (C-r) of $\Phi_r(s)$.

\                                                                               

\bigskip

\noindent
Gautami Bhowmik\\
Universit{\'e} de Lille 1, Laboratoire Paul Painlev{\'e},\\
U.M.R. CNRS 8524, 59655 Villeneuve d'Ascq Cedex, France\\
bhowmik@math.univ-lille1.fr
\bigskip

\noindent
Kohji Matsumoto\\
Graduate School of Mathematics, Nagoya University,\\
Chikusa-ku, Nagoya 464-8602, Japan\\
kohjimat@math.nagoya-u.ac.jp

\end{document}